\documentclass[letterpaper, 10 pt, conference]{ieeeconf}  

\IEEEoverridecommandlockouts                              

\usepackage{tikz}
\usepackage{pgfplots}
\tikzstyle{block} = [draw, rectangle, 
    minimum height=3em, minimum width=5em]
\tikzstyle{input} = [coordinate]
\tikzstyle{output} = [coordinate]
\tikzstyle{pinstyle} = [pin edge={to-,thin,black}]
\usetikzlibrary{shapes.geometric,shapes.arrows,decorations.pathmorphing,calc}
\usetikzlibrary{matrix,chains,scopes,positioning,arrows,fit}

\usepackage{siunitx}
\usepackage{epsfig}
\usepackage{times}
\usepackage{amssymb}
\usepackage{dsfont}
\usepackage{multirow}
\usepackage{algorithm}
\usepackage[noend]{algcompatible}

\usepackage{amsmath}
\usepackage{steinmetz}
\usepackage{subfig}
\usepackage{cite}
\usepackage{graphicx}    
\usepackage{amsfonts}
\pgfplotsset{compat=newest, ticks=none}

\title{\LARGE \bf
Control-Tutored Reinforcement Learning}

\author{F. De Lellis$^{1}$, F. Auletta$^{1}$, G. Russo$^{2}$, P. De Lellis$^{1}$ and M. di Bernardo$^{1}$
 \thanks{$^{1}$Department of EE and ICT, University of Naples Federico II, Italy {\tt\small mario.dibernardo@unina.it}}
 \thanks{$^{2}$School of EEE, University College Dublin, Ireland}
}

\begin{document}

\maketitle
\thispagestyle{empty}
\pagestyle{empty}

\begin{abstract}
We introduce a control-tutored reinforcement learning (CTRL) algorithm. The idea is  to enhance tabular learning algorithms so as to  improve the exploration of the state-space, and substantially reduce learning times by  leveraging some limited knowledge of the plant encoded into a tutoring model-based control strategy. We illustrate the benefits of our novel approach and its effectiveness by using the problem of controlling one or more agents to herd and contain within a goal region a set of target free-roving agents in the plane. 
\end{abstract}

\section{Introduction}
Reinforcement learning (RL) \cite{Sutton1998} is increasingly used to learn control policies from data \cite{kober2013reinforcement,garcia2015comprehensive,cheng2019end} in applications. While the lack of requiring a formal model of the environment/plant renders this approach particularly appealing in many applications, a key drawback is its sample inefficiency. 
Essentially, this is due to the fact that RL  finds the control policy by exploring heuristically the Markov Decision Process underlying the problem accepting possible failures. Unfortunately, long training phases are often unacceptable and failures while learning might lead to unsafe situations. 
Moreover, many applications are often  characterized by a continuous state-space and using RL requires a dense discretization of the system state space, yielding a substantial growth of learning times sometimes incompatible with the nature and scope of the control problem of interest.

Therefore, much research effort is being devoted to design safer and more sample efficient RL algorithms. An example is model-based reinforcement learning that has been used both to empower learning processes e.g. \cite{gu2016continuous, deisenroth2011pilco} and to guarantee safety in critical cases where a model is available e.g. \cite{rosolia2017learning, FerraroRusso2019}.  The introduction of some model also helps to bring some degree of stability to the overall learning process\cite{berkenkamp2017safe}. Other extensions include the Deep Q-Network (DQN) approach presented in \cite{mnih2015human} and the Actor-Critic paradigm \cite{Sutton1998},\cite{mnih2016asynchronous},\cite{lillicrap2015continuous}.

In this paper, we present an alternative model-based approach where a feedback control strategy designed with only limited or qualitative knowledge of the system dynamics is used to enhance the RL algorithm {\em when needed}.  The  resulting {\em control-tutored Q-learning} (CTQL) algorithm is better apt to deal with continuous or large state spaces while retaining many of the features of a tabular method. Our algorithm is complementary to other existing model-based approaches such as \cite{Fathinezhad2016,brunner2017repetitive}. Indeed, in our setting, the control tutor supports the process of exploring the optimization landscape by suggesting possible actions based on its partial knowledge of the system dynamics that the learning agent can take whenever it is unable to find a better action to take by itself. In so doing, the control tutor contributes to completing the Q-table speeding up the convergence of the learning process. A related  but different idea was recently presented in \cite{FerraroRusso2019}  where RL is mirrored with a Model Predictive Controller (MPC) and a different strategy is used to orchestrate transitions between RL and MPC.

To validate our approach, we apply CTQL to solve a challenging multi-agent {\em herding control problem}. The goal is that of driving  a set of agents (the herders) so that they can confine another group of autonomous roving agents (the targets) into some predefined area of the plane and keep them therein \cite{Licitra2017,pierson2017controlling}. We show that CTQL can be effectively used to solve this ``herding'' problem both in the case of one herder agent influencing one target and in the more challenging case of two herders controlling the motion of a group of ten target agents. Interestingly, we find that CTQL obtains better performance and convergence than Q-learning or feedback control on their own, solving the herding problem even when they are unable to do so.

\section{Reinforcement Learning: the key ingredients} \label{sec:RL}
Reinforcement learning is an area of machine learning which provides a set of methods that rely on approximations producing suboptimal policies \cite{Bertsekas2019} for the solution of dynamic programming problems in the presence of uncertain dynamics \cite{Sutton1998}. We briefly review here its main ingredients to properly expound the novel CTQL algorithm in the right context.


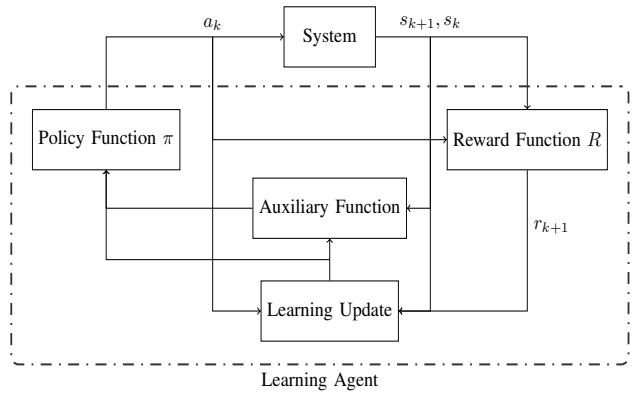
\begin{figure}[tbp] 
\centering
\resizebox{8.4cm}{5.2cm}{
\begin{tikzpicture}[auto, node distance=1.8cm]
    \node [input, name=input] {};
    \node [block, right of=input] (system) {{System}};
    \node [block, below of=system, right of=system, xshift=2cm] (RF) {Reward Function $R$};
    \node [block, below of=system, yshift=-1.2cm] (AF) {Auxiliary Function};
    \node [block, below of=AF] (LU) {Learning Update};
    \node [block, below of=system, left of=system, xshift=-2.5cm] (policy) {Policy Function $\pi$};

    \draw [->] (system) -| node [pos=0.18] [name=S] {$s_{k+1},s_{k}$} (RF);
    \draw [->] (RF) |-  node [pos=0.2] {$r_{k+1}$} (LU);
    \draw [->] (S) |- (AF);
    \draw [->] (AF) -| (policy);
    \draw [->] (policy) |- node [pos=0.8] [name=A] {$a_{k}$} (system);
    \draw [->] (A) |- (RF);
    \draw [->] (S) |- (LU);
    \draw [->] (A) |- (LU);
    \draw [->] (LU) -- node [pos=0.5] [name=temp] [right] {} (AF);
    \draw [->] (temp) -| (policy);
    
    \tikzset{black dotted/.style={draw=black!80!white, line width=1pt, dash pattern=on 1pt off 4pt on 6pt off 4pt, inner sep=4mm, rectangle, rounded corners}};

    \node (first dotted box) [black dotted, fit = (RF) (policy) (AF) (LU)] {};
    \node at (first dotted box.south) [below] {Learning Agent};
\end{tikzpicture}
}
\caption{Reinforcement Learning scheme} 
\label{fig:RL}
\end{figure}
%
A RL control loop is schematically shown in Fig. \ref{fig:RL}. Specifically, the interaction between the control agent and the environment/plant (or simply system in what follows) is described through: (i) the state space $\mathcal{S}$ containing all possible system states; (ii)  the action space $\mathcal{A}$ of all possible actions the agent can take to influence the system state. As shown in Fig.~\ref{fig:RL}, other key components of a RL control algorithm are: (i)  a policy function $\pi$;  (ii) a reward function $R$; (iii) a learning update rule; and (iv) an auxiliary function. 

The policy selection function $\pi$ is used to determine what action $a\in \mathcal{A}$ to apply to the system starting from state $s\in\mathcal{S}$ at time $k$. The effects of such action, say $a_k$, are evaluated via the reward function that evaluates the effects of that action with respect to the control goal. Namely,  given the action taken at step $k$, the current state of the system $s_k$, and the computed next state value $s_{k+1}$, the expected reward $r_{k+1}$ is computed. The learning update rule is then used to update an auxiliary function storing the expected rewards for taking a certain action when the system is in a given state. Such an auxiliary function is interrogated at each step by the learning agent to decide what action to take next and represents the ``experience" accumulated by the agent during the learning process.

In the Q-learning approach \cite{Sutton1998}, \cite{Watkins1992}, the {\em auxiliary function} is expressed as a tabular function, $Q(s,a)$ with $a \in \mathcal{A}$, $s \in \mathcal{S}$, called Q-table.
The state and action spaces are assumed to be discrete and of finite cardinality.
In this way, the learning agent accesses the table by the state, at time $k$, and selects the best action to take according to the values stored in the Q-table. The {\em policy selection function} exploits the $\varepsilon$-greedy criterion \cite{Sutton1998} and is defined as follows:
\begin{equation} \label{eq:piQ}
    \pi_Q(s_k) =     \begin{cases} 
        \arg \max\limits_{a \in \mathcal{A}} \{{Q}(s_k,a)\},  &   \text{with probability ($1 -\varepsilon$)}\\
        \mbox{rand}\ (a),   &   \text{with probability $\varepsilon$}
    \end{cases}
\end{equation}
 where $\varepsilon$ is a positive constant in the range $]0,1[$ representing the probability of taking a random action instead of an action stored in the Q-table. Randomness in the policy promotes exploration and fulfills the hypotheses needed to prove convergence of the algorithm towards the optimal solution \cite{Watkins1992}. The \emph{learning update rule} is defined as follows:
 %
\begin{equation} \label{eq:LUP}
\begin{split}
    Q_{k+1}(s_{k},a_{k})= &\ Q_k(s_{k},a_{k}) + \alpha[r_{k+1} + \\
                          & + \gamma \max_{a \in \mathcal{A}}Q_k(s_{k+1},a) - Q_k(s_{k},a_{k})]
\end{split}
\end{equation}
where $r_{k+1}$ is the reward obtained by selecting action $a_{k}$ from state $s_{k}$, and $s_{k+1}$ is the next state, and $Q_k$ is the value stored in the Q-table at time $k$.
The parameters $\alpha$ and $\gamma$ are both in the range $[0,1]$ and are known as the learning rate and the discount factor, respectively.

As mentioned in the introduction, the main problem of this algorithm is related to the need of discretizing the action and state spaces that,  for continuous dynamical systems, can lead to a substantial growth of the Q-table and poor learning performance \cite{lillicrap2015continuous}. Hence, the need for enhanced RL algorithms as the one we propose next.



\section{Control Tutored Q-learning (CTQL)} \label{sec:CTQL}
\begin{figure}[t]
\centering
\resizebox{8.4cm}{6cm}{
\begin{tikzpicture}[auto, node distance=1.8cm]
    
    \node [input, name=input] {};
    \node [block, right of=input] (system) {System};
    \node [block, below of=system, right of=system, xshift=1cm] (RF) {Reward Function R};
    \node [block, below of=system, yshift=-2cm, text width=3.1cm] (Q) {Auxiliary Function: Q-Table};
    \node [block, below of=Q] (LU) {Learning Update};
    \node [block, below of=system, left of=system, xshift=-2cm] (policyQ) {\large $\pi_Q$};
    \node [block, below of=system, left of=system, xshift=-4.4cm] (policyT) {\large $\pi_T$};
    \node [block, below of=policyT, yshift=-0.4cm] (CL) {Control Law};

    \draw [->] (system) -| node [pos=0.25] [name=S] {$s_{k+1},s_{k}$} (RF);
    \draw [-] (system) -- ++(4.7cm,0) coordinate(S1) {};
    \draw [->] (RF) |-  node [pos=0.1] {$r_{k+1}$} (LU);
    \draw [->] (S1) |- (Q);
    \draw [->] (Q) -| (policyQ);
    \draw [<-o] (system.west) -- node [pos=0.64] [name=A] [above] {$a_{k}$} ++(-1.6cm,0) coordinate(p1){};
    \draw [->] (A) |- (RF);
    \draw [->] (S1) |- (LU);
    \draw [->] (A) |- (LU);
    \draw [->] (LU) -- node [pos=0.5] [name=temp] [right] {} (Q);
    \draw [->] (temp) -| (policyQ);
    \draw [->] (S1) |- ([xshift=-0.3cm,yshift=-0.27cm]LU.south west) -| (CL);
    \draw [->] (CL) -- node [pos=0.5] [name=v] [right] {{$v_k$}} (policyT);
    \draw [-] (policyQ.north) -- ++(0,1cm) coordinate(pq1){};
    \draw [-o] (pq1) -- ++(0.7cm,0) coordinate(pq2){};
    \draw [-] (policyT.north) -- ++(0,1.5cm) coordinate(pt1){};
    \draw [-o] (pt1) -- ++(3.08cm,0) coordinate(pt2){};
    \draw [-] (p1) -- (pq2);
    \draw[<-] ($(pt2)+(0.3cm,0.1cm)$) to [bend right]($(pq2)+(0.3cm,-0.1cm)$);
    \tikzset{black dotted/.style={draw=black!80!white, line width=1pt, dash pattern=on 1pt off 4pt on 6pt off 4pt, inner sep=5mm, rounded corners}};
    \tikzset{blue dotted/.style={draw=blue!80!white, line width=1pt, dash pattern=on 1pt off 4pt on 6pt off 4pt, inner sep=1.5mm, rectangle, rounded corners}};
    \tikzset{red dotted/.style={draw=red!80!white, line width=1pt, dash pattern=on 1pt off 4pt on 6pt off 4pt, inner sep=1.1mm, rounded corners}};
    \tikzset{green dotted/.style={draw=green!70!black, line width=1pt, dash pattern=on 1pt off 4pt on 6pt off 4pt, inner sep=1.8mm, rounded corners}};

    \draw[-] [black dotted] ($(RF)+(2.1cm,1cm)$) |- ($(LU)+(0.1cm,-1cm)$);
    \draw[-] [black dotted] ($(LU)+(0.1cm,-1cm)$) -| node [pos=0.05, yshift=0.25cm] [name=LA] {Learning Agent} ($(CL)+(-1.4cm,0.1cm)$);
    \draw[-] [black dotted] ($(CL)+(-1.4cm,0.1cm)$) -- ($(policyT)+(-1.4cm,0.2cm)$) |- ($(p1)+(0cm,0.6cm)$) -- ++(0.3cm,0cm) coordinate(x1) {};
    \draw[-] [black dotted] (x1) |- ($(RF)+(2.1cm,1cm)$) {}; 
    
    \node (second dotted box) [blue dotted, fit = (CL)] {};
    \node at (second dotted box.south) [below, xshift=-0.7cm]  {Tutor};
    
    \node (third dotted box) [red dotted, fit = (Q) (LU) (RF)] {};
    \node at (third dotted box.south) [above, xshift=1.2cm]  {Q-learning};
    
    \node (fourth dotted box) [green dotted, fit = (policyT) (policyT) (p1) (pt2)] {};
    \node at (fourth dotted box.south) [above, xshift=0.06cm, yshift=2.25cm]  {Policy Selection $\pi$};
\end{tikzpicture}
}
\caption{Control Tutored Q-learning (CTQL) Schematic} 
\label{fig:CTQL}
\end{figure}
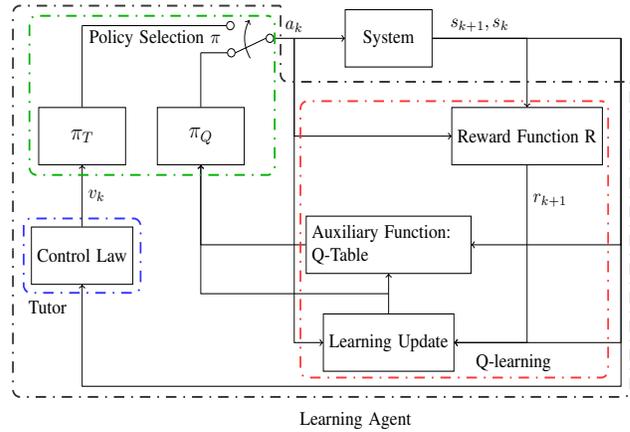

The key idea behind CTQL is schematically summarized  in Fig.~\ref{fig:CTQL}. 
Specifically, at each time step $k$, the learning agent selects its next action $a_k$ from a given system state $s_k$,
by choosing either the action suggested by the control tutor via a model-based policy $\pi_T$, or the one suggested by the standard $\varepsilon$-greedy policy $\pi_Q$ used for the Q-learning as defined  in \eqref{eq:piQ}.  In so doing, CTQL adopts the same Q-table structure and learning update rule of Q-learning,  but exploits a new policy selection function, say $\pi$.

Mathematically, the policy selection function in CTQL is a switching policy defined as:
\begin{equation} \label{eq:4-6}
    \pi(s_k) = 
    \begin{cases} 
        \pi_Q(s_k),   &    \max\limits_{a \in \mathcal{A}}\{Q(s_k,a)\}>0,\\
        \pi_T(s_k),   &   \text{otherwise}.
    \end{cases}
\end{equation}

According to this policy, at step $k$, given the state $s_k$, the learning agent checks the entries of the Q-table for all actions $a \in \mathcal{A}$. If at least one of these entries is positive, then that action $a_k$ is selected by Q-learning, otherwise the action is chosen that is suggested by the control-tutor via the policy  $\pi_T(s_k)$. This policy is defined as follows:
\begin{equation} \label{eq:4-5}
  \pi_T(s_k) = 
    \begin{cases} 
        \arg \min\limits_{ a \in \mathcal{A}}\{\|v_k - a\|\}, &   \text{with probability ($1 -\varepsilon$)},\\
        \mathrm{rand}(a),   &   \text{with probability $\varepsilon$},
    \end{cases}
\end{equation}
where $v_k$ is the control input  generated by the control tutor using a feedback controller designed on a rough model of the plant. As such input does not necessarily belong to $\mathcal A$, the policy function $\pi_T$ selects the action $a\in\mathcal{A}$ which is closest to $v_k$.

%

Once, the action is selected from either $\pi_Q(\cdot)$ or $\pi_T(\cdot)$, the corresponding expected reward is then computed and used to update the Q-table. The pseudocode of the CTQL algorithm is given in Algorithm \ref{alg:CTQL}. Note that both $\pi_Q(\cdot)$ and $\pi_T(\cdot)$ contain some degree of randomness to guarantee that when implemented the policy selection function of the CTQL is still within the scope of the probabilistic proof of convergence available for the Q-learning algorithm and described in \cite{Watkins1992}.


\begin{algorithm}[htbp]
  \caption{Control Tutored Q-learning}
  \label{alg:CTQL}
  \begin{algorithmic}
    \STATE {Initialize $Q(s,a) = 0,\forall s \in \mathcal{S},a \in \mathcal{A}$}
    \FOR{$n = 1$ to $N_{tr}$}
      \STATE{Detect intial state $s_0$}
      \FOR{$k = 0$ to $T_{tr}$}
        \IF{$\max\limits_{a \in \mathcal{A}}\{Q(s_k,a)\}>0$}
           \STATE{$a_k \gets\pi_Q(s_k)$}
       \ELSE
           \STATE{$a_k \gets\pi_T(s_k)$}
       \ENDIF
       \STATE {Observe $r_{k+1}$ and $s_{k+1}$}
       \STATE {$Q(s_{k},a_{k})\gets(1-\alpha)Q(s_{k},a_{k}) + \alpha[r_{k+1} + $}
       \STATE {$\ \ \ \ \ \ \ \ \ \ \ \ \ \ \ \ \ \ \ \ + \gamma \max\limits_{a \in \mathcal{A}}Q(s_{k+1},a)]$}
      \ENDFOR
    \ENDFOR
  \end{algorithmic}
\end{algorithm}

To illustrate the viability and effectiveness of CTQL, we apply it to solve the herding problem in robotics and discuss its performance by comparing it to a traditional (untutored) Q-learning approach.
\section{Application to the Herding Problem}
We consider the problem of letting one or more mobile agents in the plane (the herders) drive the motion of a group of autonomous agents (the targets) so as to move them towards some goal region and mantain them therein.
Under the assumption that the herders only possess limited knowledge of the dynamics of the targets, we will solve the problem of controlling the herders by using CTQL. For the sake of simplicity, we assume all agents are able to adjust their velocities almost instantaneously, as done for example in \cite{albi2016invisible}. In what follows we will use the pedix `$\tau$' to denote quantities pertaining to the target agents and the pedix `$h$' for those concerning the herders.



\subsection{Problem Formulation}
Assuming, the target agents' velocity is upper bounded by some maximum velocity $v_{\tau}^{\max}$, the dynamics of the target agents is assumed to be: 
\begin{equation} \label{eq:3-3}
\dot{x}^i_\tau(t) =
    \begin{cases}
         f_i(x_\tau^i(t), x_h(t)),     &    \|f_i(x_\tau^i(t), x_h(t))\| < v_{\tau}^{\max}\\
        v_{\tau}^{\max}e^{\jmath \scalebox{0.74}{\phase{f_i(x_\tau^i(t), x_h(t))}}}, &   \text{otherwise}
    \end{cases}
\end{equation}
where $\jmath$ is the imaginary unit, $ x_\tau^i(t) \in \mathbb{R}^2$ is the position of the $i$-th target agent (out of $N$) at time $t$, $x_h(t) = [x_h^1(t),\ldots, x_h^M(t)]^T\in\mathbb{R}^{2M}$ is the vector stacking the positions  of the $M$ herder agents at time $t$, and the vector field $f_i:\mathbb{R}^{2(M+1)}\mapsto\mathbb{R}^2$ is the sum of two contributions, i.e. $f_i=f_1+f_2^i$.

Here, the term $f_1$ models the action of the herders onto the target and is assumed to be the same for all the targets. It is defined as:
\begin{equation}
     f_1(x_\tau^i, x_h)  := \beta_1\sum_{j=1}^M \frac{ x_\tau^i -  x_h^j}{\| x_\tau^i -  x_h^j\|^3} U ( x_\tau^i,  x_h^j, \rho_\tau)
\end{equation}
where we omitted the explicit dependence on time $t$, $\rho_\tau$ is the targets' influence radius, $\beta_1>0$ is a constant gain modelling the intensity of the coupling with the herder, and $U$ is an interaction function defined as
\begin{equation}\label{eq:step}
    U ( x_\tau^i,  x_h^j, \rho_\tau) = \begin{cases}
    1, & \| x_\tau^i -  x_h^j\| < \rho_\tau \\ 
    0, & \text{otherwise} 
    \end{cases}
\end{equation}
that ensures that the coupling between target and herder agents is active only if their relative distance is smaller than some $\rho_\tau>0$.

The term $f_2^i$  represents the target own random dynamics defined as:
\begin{equation} 
      f_2^i := \beta_2^i(t) e^{\jmath \theta^i(t)}
\end{equation}
where $\beta_2^i(t)$ and $\theta^i(t)$ are scalars updated every $\Delta t$ seconds with values extracted from uniform distributions $ \mathcal{U}(0,\beta_{\max}) $ and $\mathcal{U}(0,2\pi)$, respectively.

The herders' speed, as for the targets, is saturated to a maximum fixed value $v_{h}^{\max}$ so that their dynamics can be written as:
\begin{equation} \label{eq:3-4}
\dot{ x}_h^j(t) = 
    \begin{cases}
         u^j(t),   &  \| u^j(t)\| < v_{h}^{\max}\\
        v_{h}^{\max}e^{\jmath \scalebox{0.8}{\phase{ u^j(t)}}}, &   \text{otherwise}
    \end{cases}
\end{equation}
where $u^j(t)$ is a control input at time $t$ to be determined in order to fulfill the control goal.

 The control objective is to design the input vector $u = [  u^1,..., {u^M}]^T$ able to drive the targets to reach and remain in the circular goal region $ G := \{  x\in {\mathbb{R}}^2 : \| x -  x_g\|<\rho_g\}$  of center $ x_g$  and radius $\rho_g$, that is, to guarantee that
\begin{equation}
    \limsup_{t\to+\infty} \| x_\tau^i(t) -  x_g\| < \rho_g, \quad \text{for all } i = 1,...,N.
\end{equation}

\subsection{Control Design}
 For the sake of simplicity, we start by considering the case where $N=M=1$ (dropping the suffixes $i$ and $j$) and the goal region is centered at the origin, i.e. $x_g=0$. We suppose the herder knows the position of the target but possesses only a conservative estimate, $\hat \rho_\tau<\rho_\tau$ of the target's true influence radius $\rho_\tau$.

We design the control input $u$ driving the herder as follows. At time $t$, if $\|x_\tau(t)-x_h(t)\|>\hat \rho_\tau$, then the herder moves towards the target at its maximum speed to reduce its distance until entering the estimated influence region at some time $\tilde t$ when $\|x_\tau(\tilde t)-x_h(\tilde t)\|\leq\hat \rho_\tau$.  Within this region the herder adopts a learning strategy to push the herder towards the goal region. For the sake of comparison, we first test how Q-learning performs to solve the problem and then move to CTQL.

\subsubsection{Q-learning implementation} \label{sec:QHerding}
We start by applying the classical Q-learning algorithm with the following definitions of state and action spaces, and of the reward function.

The action and state space are defined as follows.  $\mathcal S:= \mathcal{D}\times\mathcal{W}\times\mathcal{V}$, where (i) $\mathcal{D}$ is the set of distances of the herder and the target from the center of the goal region; (ii) $\mathcal{W}$ is the set of angular positions of the herder; (iii) $\mathcal{V}$ is the set of possible speeds of the  target. In our implementation the sets $\mathcal{D},\mathcal{W}, \mathcal{V}$ are discrete sets and are defined in detail in the Appendix. The action space (see also the Appendix) is the set of possible discretized values of the input vector $u$ to the herder dynamics given by \eqref{eq:3-4}.



Let $ x_{\tau,k}$ $(x_{h,k})$ be the position of a generic target (herder) agent at a discrete  time instant $k$. 
Then, the reward function implemented in our experiments is: 
\begin{equation}
    R(a_k,s_k,s_{k+1}) := k_1R_1(a_k,s_k,s_{k+1}) + k_2R_2(a_k,s_k,s_{k+1})
\end{equation}
where
\begin{equation}
R_1(a_k,s_k,s_{k+1}) =  (| x_{\tau,k}| - | x_{\tau,k+1})|),
\end{equation}
\begin{equation}
R_2(a_k,s_k,s_{k+1}) = \sigma(\Bar k(| x_{h,k+1}| - \rho_g)) - 1,
\end{equation}
with $k_1$, $k_2$, and $\Bar k$ being positive constant gains, and where $\sigma(z) := {1}/{(1 - e^{-z})}$ was chosen w.l.o.g. as a decreasing function of its argument.

\begin{figure}[htbp]
\begin{center}
\includegraphics[width=10cm]{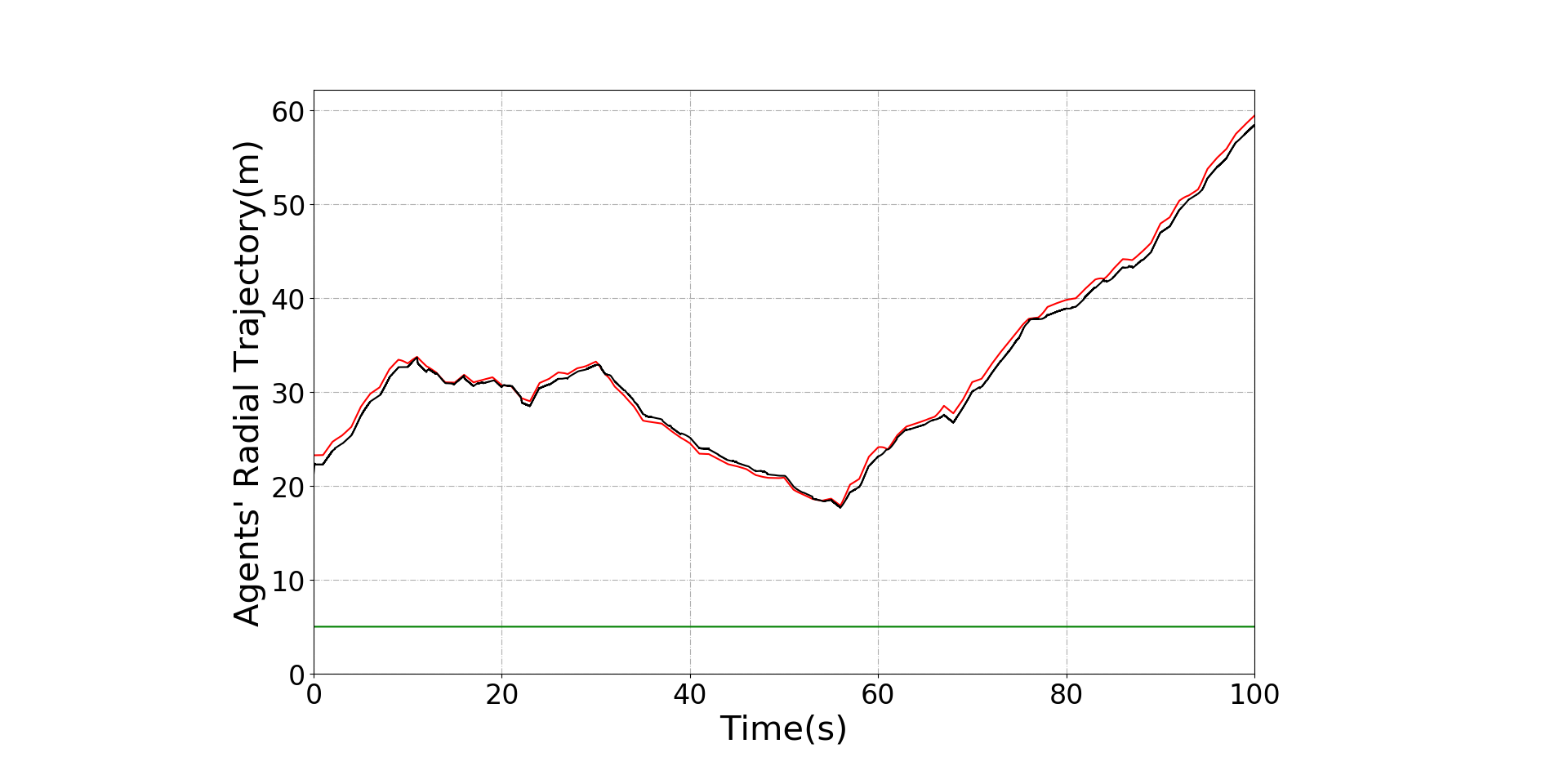} 
\includegraphics[width=10cm]{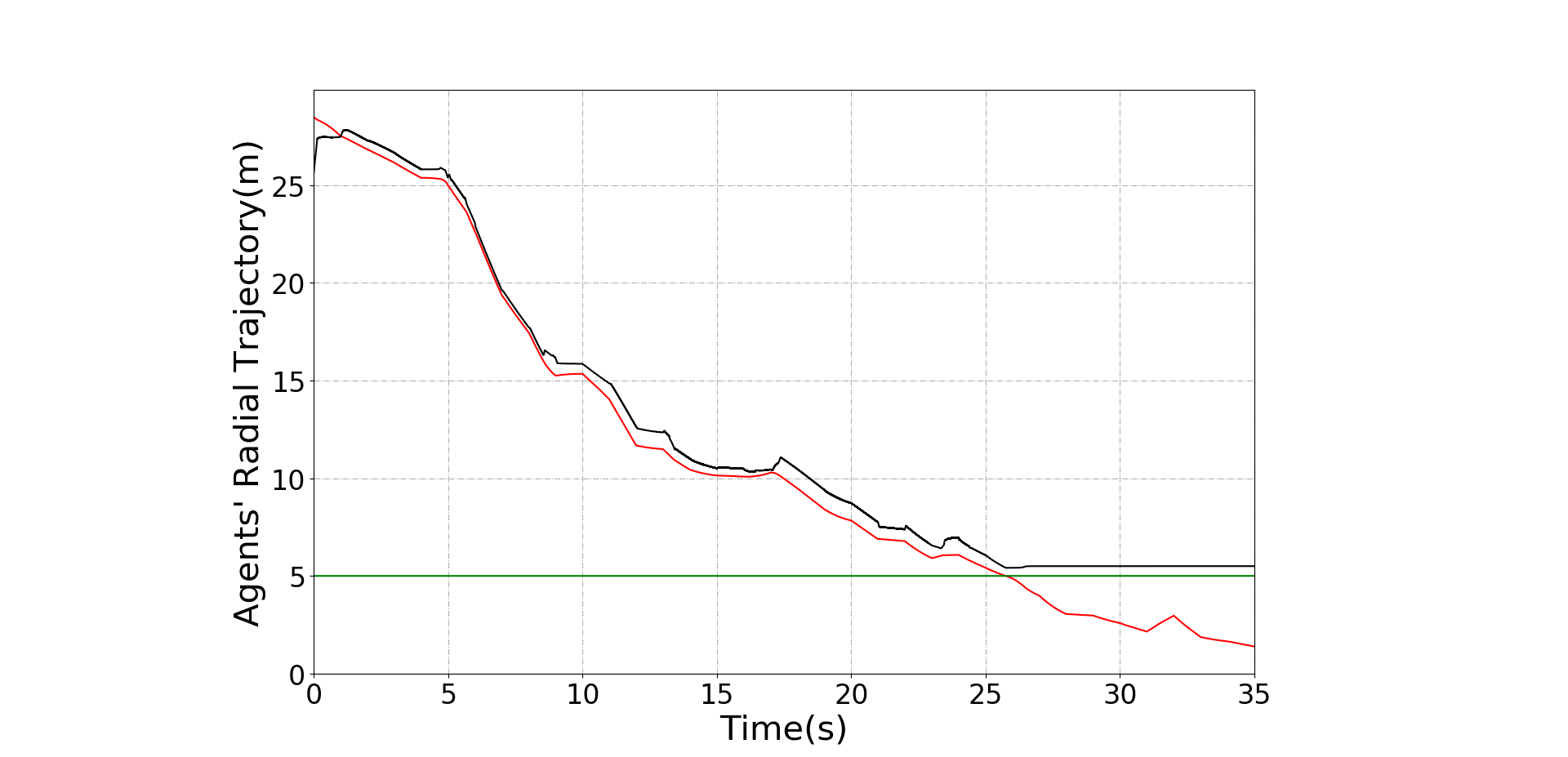}
\caption{Performance of the Q-learning algorithm after $N_{tr}:=5000$ trials and the state discretization is finer (top panel) or coarser (bottom panel). The radial distance of the herder (black line) and the target (red line) are shown together with the radius of the goal region (green line).} \label{fig:AQT}
\end{center}
\end{figure}

We test the Q-learning algorithm by considering two discretizations of state space, a finer and a coarser one (see the Appendix for further details). When using the finer discretization, we see that, as shown in Fig.~\ref{fig:AQT} (top panel), Q-learning is unable to achieve the control goal after over 5000 training trials. Convergence is instead achieved within about $25\si{\second}$ after a training phase of about the same duration when a coarser discretization of the state space is used, see Fig. \ref{fig:AQT} (bottom panel). 

\subsubsection{Control-Tutored Q-learning Application} 
We  move next to adopting the CTQL approach described in Sec.~\ref{sec:CTQL}. The state and action spaces, as well as the reward function, are the same as those proposed in Sec.~\ref{sec:QHerding}.

The design of the tutoring control law requires some model of the expected dynamics of the targets. 
We assume that only an estimate of the target true dynamics is available which we suppose to be given by the inaccurate model:
\begin{equation}\label{eq:4-3}
\dot{ x}_\tau(t) = \gamma( x_\tau(t) -  x_h(t)) U( x_\tau(t), x_h(t),\hat \rho_\tau),
\end{equation} 
where $\gamma>0$ is a gain modelling the intensity of the coupling between the target and the herder, $\hat \rho_\tau \leq \rho_\tau$, and $U(\cdot)$ is the step function defined in~\eqref{eq:step}.

Assuming the target's dynamics as in \eqref{eq:4-3}, we then select the herder control input so as to push the target position towards the origin; namely we choose
\begin{gather} \label{eq:4-7}
u(t) = k_i\dot{x}_\tau(t) + k_p x_\tau(t),
\end{gather}
where $k_i$ and $k_p$ are two positive control gains.

With this choice of $u$ when the target and the herder interact, the target dynamics becomes:
\begin{gather} \label{eq:4-8}
    \ddot{x}_\tau(t) = \dot{x}_\tau(t) - \dot{x}_h(t) = (1-k_i)\dot{x}_\tau(t) - k_p x_\tau(t),
\end{gather}
so that any choice of $k_i>1$ and $k_2>0$ would achieve convergence to the origin were the dynamics \eqref{eq:4-3} the correct ones. Without loss of generality here we choose $k_i = 2$, $k_p = 0.1$.

\begin{figure}[htbp]
\begin{center}
\includegraphics[width=10cm]{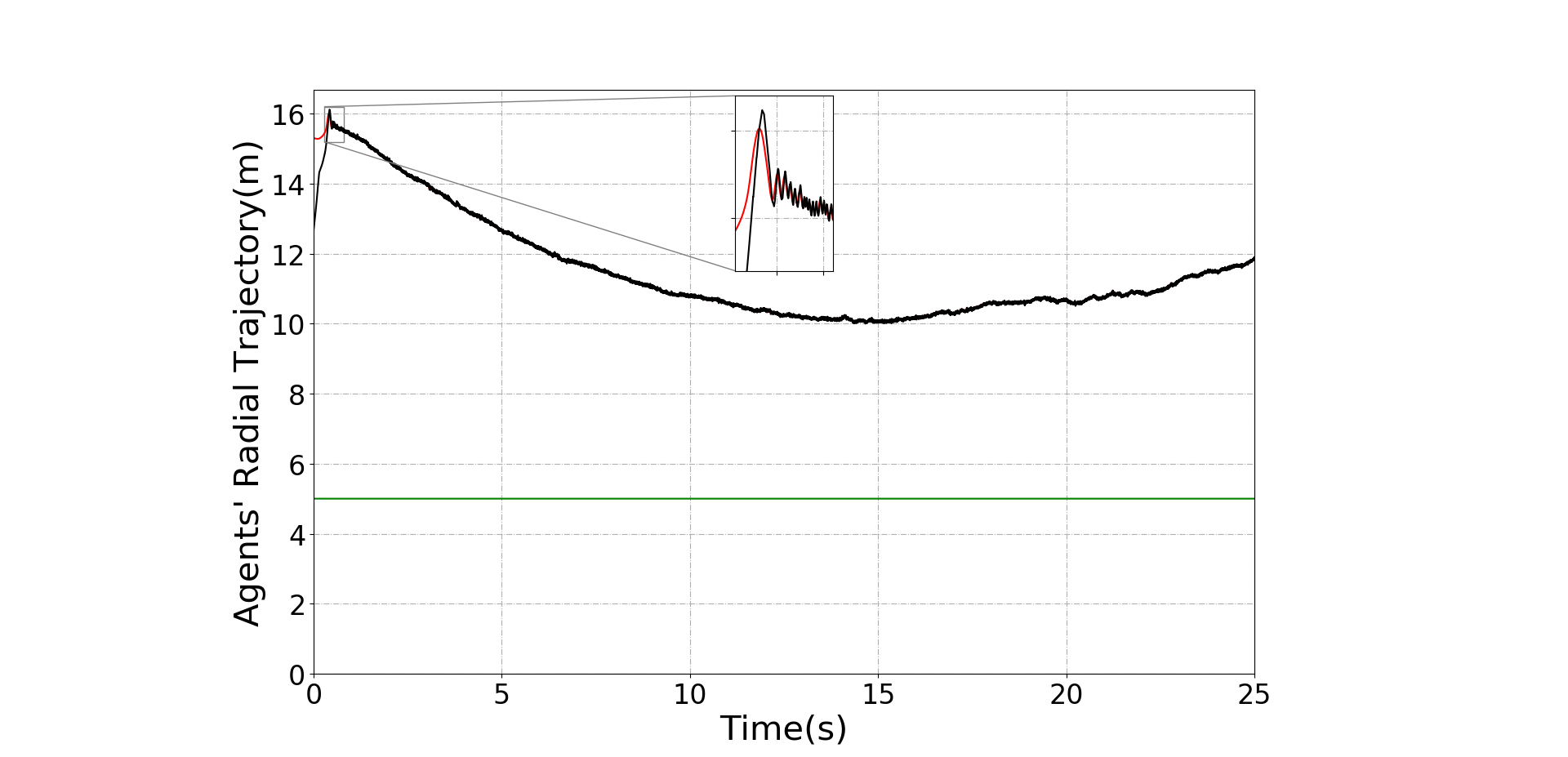}    
\caption{Performance of the control tutor without any learning. The inset shows a zoom of the transient dynamics during the interval $t \in [0.28, 0.8]\si{\second}$. The color codes are described in the caption of Fig. \ref{fig:AQT}.} 
\label{fig:AT}
\end{center}
\end{figure}

As expected, when applied to control the ``true" target dynamics, we observe that, as shown in Fig. \ref{fig:AT}, the herder driven by \eqref{eq:4-7} fails to achieve the desired goal as the target escapes the region where they actually interact and becomes uncontrollable. (Note that a better choice of the controller or the gains might resolve this issue for the approximate model but here leave the controller unchanged as we wish to explore whether our CTQL approach can instead achieve convergence even when the control tutor is designed on a set of very simplifying qualitative assumptions such as those we made.)
\subsection{Numerical Validation}
\subsubsection{CTQL herding of a single target}
\begin{figure}[tb]
\begin{center}
\includegraphics[width=10cm]{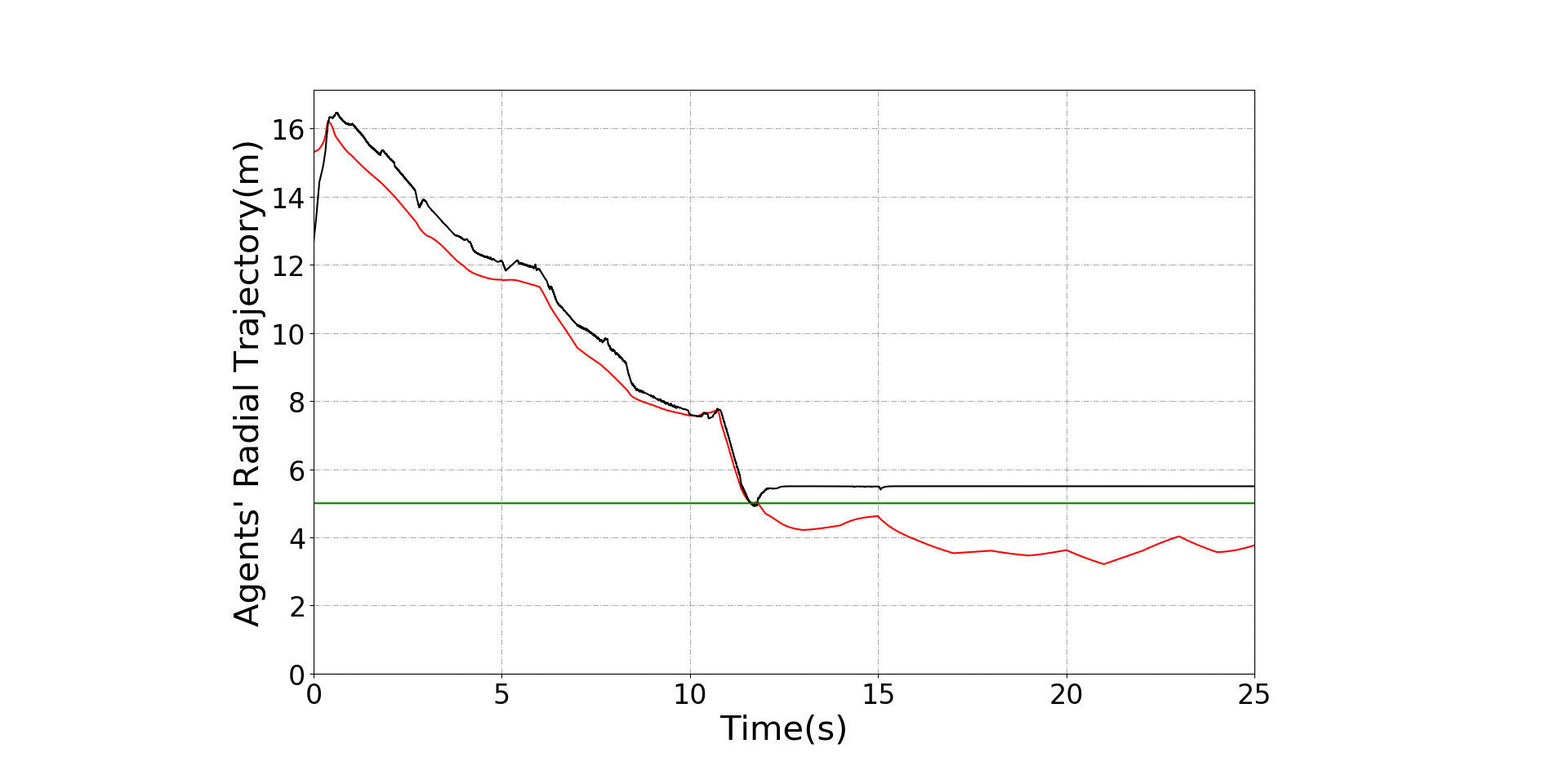} 
\includegraphics[width=10cm]{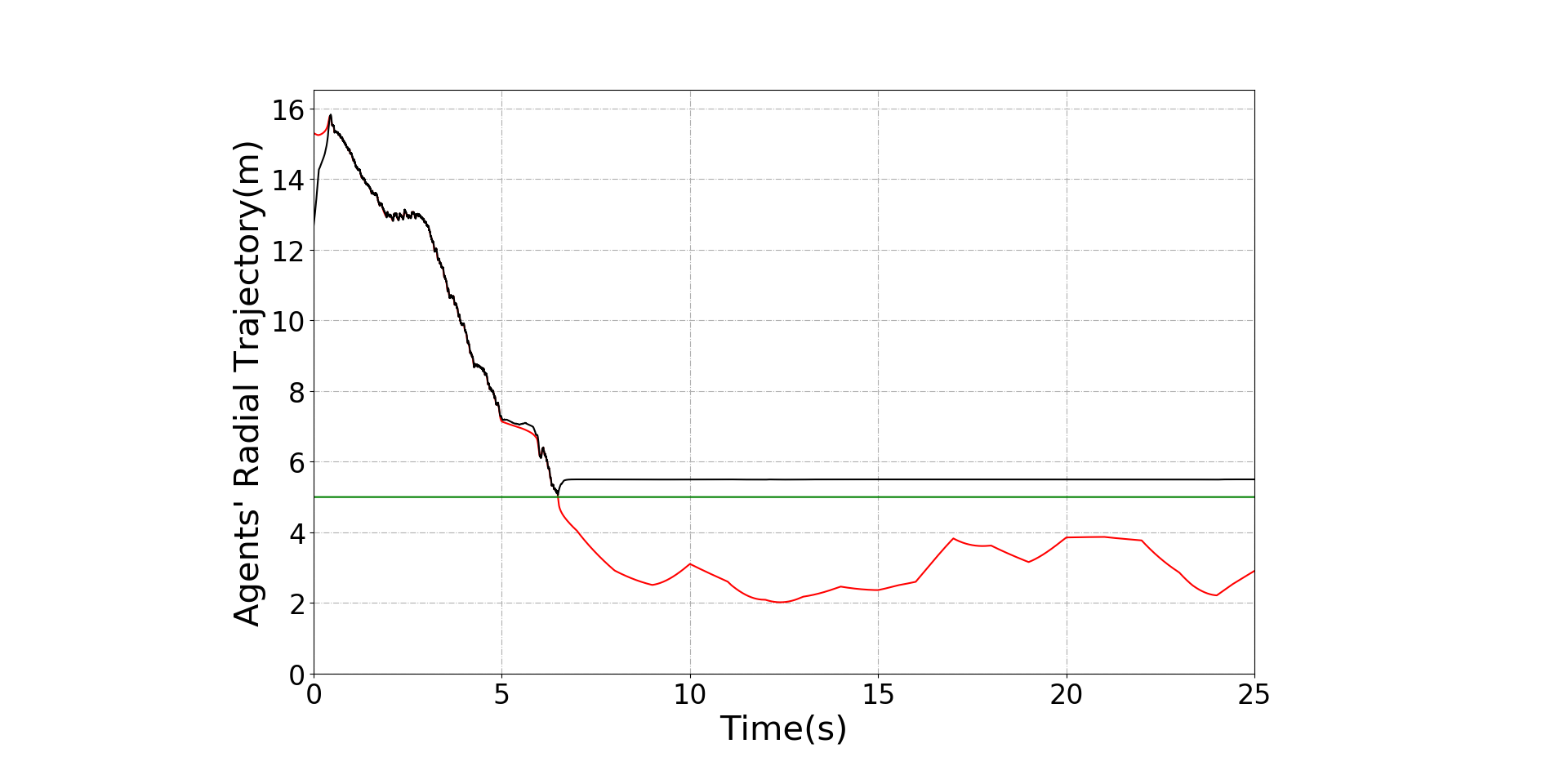}     
\caption{Performance of the CTQL algorithm with a finer (top panel) or coarser (bottom panel) state space discretization after just 1 learning trial. Convergence is immediately achieved in both cases. The color codes are described in the caption of Fig. \ref{fig:AQT}.} 
\label{fig:ATTVFB}
\end{center}
\end{figure}

As shown in Fig. \ref{fig:ATTVFB}, in the case of one herder interacting with one target, using the CTQL approach with the control tutor designed above is successful after just one or two training trials independently of the state discretization used. A summary of the performance and convergence times of CTQL compared with those where the Control Tutor (CT) or Q-learning (QL) are used on their own is shown in Table \ref{tab:narrow}. The numerical experiments where initiated with random initial conditions, $x_\tau(0)$ uniformly selected in $[15,30]$, and $x_h(0)$ such that the initial distance $\|x_\tau(0) - x_h(0)\|$ is uniformly distributed in $[\hat \rho_\tau, \hat \rho_\tau + 2]$.
We observe that the control tutor on its own (without learning) is never successful while the Q-learning performance strongly depends on the state discretization used. CTQL  instead always achieves convergence guaranteeing robustness to hyperparameters selection such as  state discretization and a very limited number of learning trials (1 or 2 in the case we tested as compared to over 5000 for Q learning with a coarse state discretization).

\begin{table}[tbp]
\caption{Performance comparison between CT, QL and CTQL with the finer (and coarser) state discretization}
\label{tab:narrow}
\begin{center}
\begin{tabular}{||c||c||c||}
\hline
Control type & success rate over 100 trials & $<t_s>$\\
\hline
QL after 5000 trials & 0\% (100\%) & $+\infty$ (37.2s) \\
\hline
CT & 0\% (0\%) & $+\infty$ ($+\infty$)  \\
\hline
CTQL & 100\% (100\%) & 14.3s (35.7s)\\
\hline
\end{tabular}
\end{center}
\end{table}

\subsubsection{CTQL of multiple herders and targets}
To further test our strategy, we considered the case of $M$ herders controlling $N>M\geq1$  targets.
In this context, herders' behavior needs to include some cooperation rule to successfully drive and contain the targets. Here, the herders use CTQL and cooperate to fill in the same Q-table. 

We assume  each herder is always aware of the current positions of all the targets. Then, using this information, herders (i) compute the center of mass (CoM) of the positions of the targets and  (ii)  split the plane into $M$ circular sectors centered at the origin by starting with the line passing through the origin and the computed CoM. Each herder then assumes control of one of such sectors taking the task of searching and recovering targets that are located in that area. The sectors are re-computed and re-allocated every $10\si{second}$.
Such division of the region of interest forces each herder to choose the targets to chase only in its sector of competence and, consequently, avoids interference among herders. 

As the velocities of all the agents are comparable, herders may end up  continually switching between two or more targets to chase without pushing any of them towards the goal region. To avoid such a case, the following rule has been introduced for herder agents: 
\begin{enumerate}
  \item Select the furthest target $\tau$ from the goal region  in your sector of competence;
  \item while trying to contain $\tau$ in the goal region $G$, check if another target, say $\tau'$, becomes the new furthest target  from $G$. If such a target exists, then
  \item compute the distance $\| x_{\tau,\tau'}\|$ between $\tau$ and $\tau'$, and
  \item if $\| x_\tau\|>\| x_{\tau,\tau'}\|$ switch the control law to contain $\tau'$ in $G$, otherwise keep containing target $\tau$.
\end{enumerate}

\begin{figure}[htbp]
\includegraphics[width=10cm]{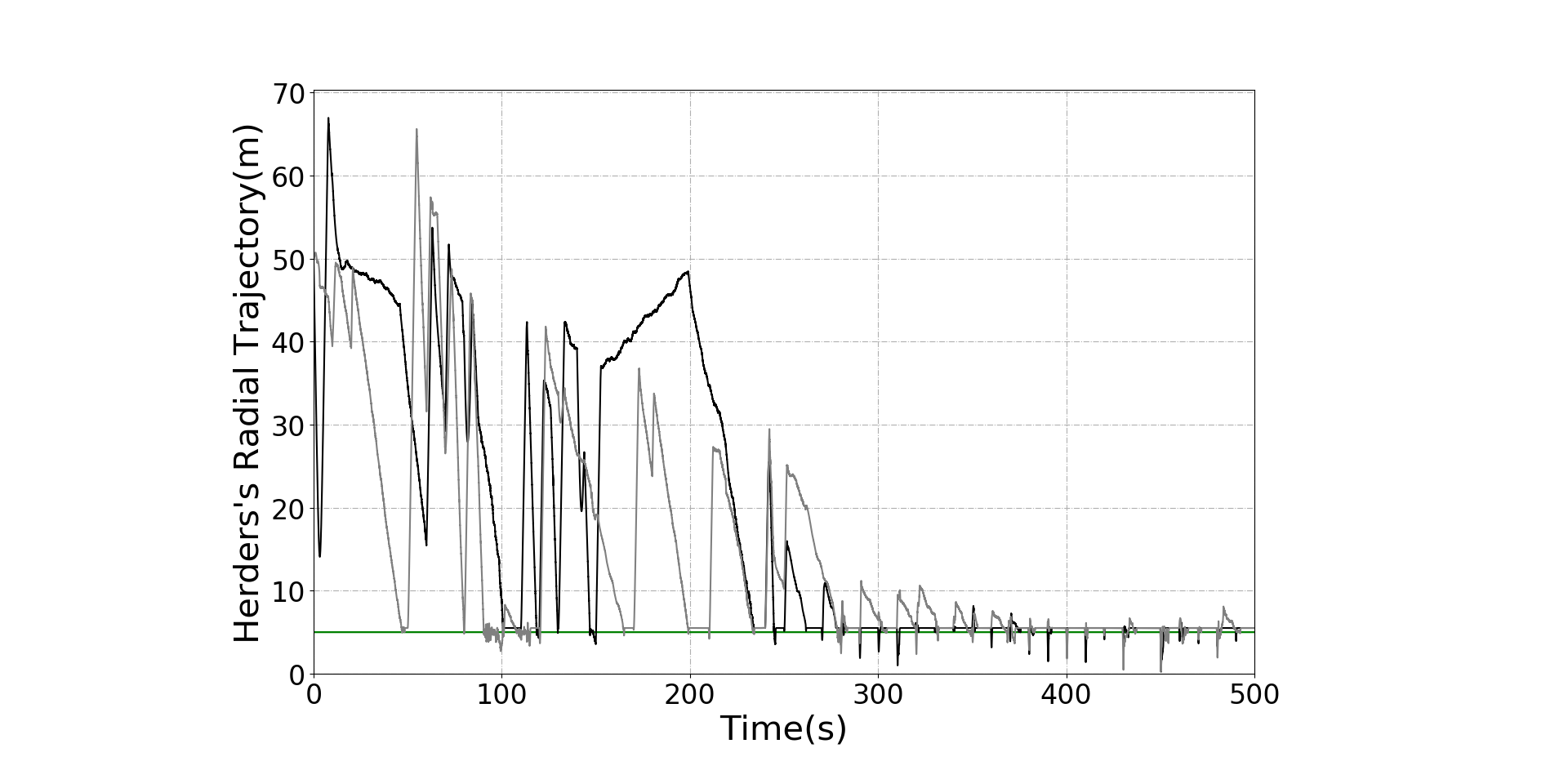}\\
\includegraphics[width=10cm]{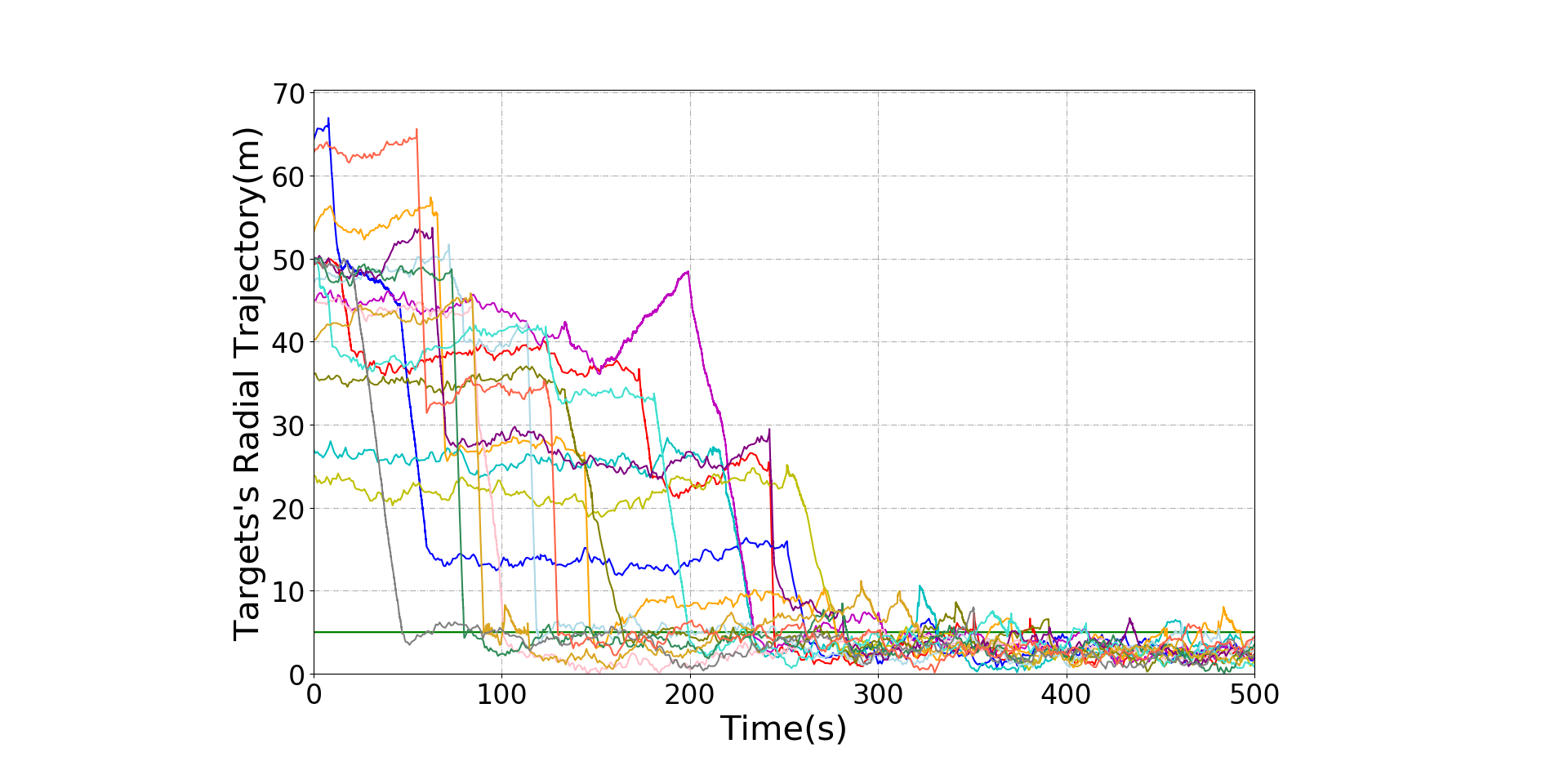}
\caption{Radial coordinates of (a) the two herding agents and (b) the 15 target agents after just one learning trial when CTQL is employed to drive the herders. The radius of the goal region is shown as a green solid line.}
\label{fig:ATTVFB10}
\end{figure}

Figure \ref{fig:ATTVFB10} shows the performance of CTQL when  $M=2$ herders interact with a group of $N = 15$ targets confirming the effectiveness of using a control tutor that allows the learning algorithm to achieve convergence after just one learning trial.

\section{Conclusions}
In this paper, we introduced an extension of Q-learning where the policy selection function is enhanced by means of a control tutor that, using a feedback control law with limited knowledge of the system dynamics, is able to support the exploration of the optimization landscape guaranteeing better convergence and shorter learning times. To illustrate the effectiveness of the approach, we discussed its application to the herding problem showing that the combination of learning and feedback control can achieve ambitious control goals even in those cases where neither would work on its own. We envisage that a similar control tutored approach can be used to enhance the performance and convergence of other more sophisticated learning algorithms. We wish to emphasize that from a control viewpoint, the combined presence of RL and feedback control renders viable the use of a control strategy that would otherwise be useless without the presence of learning. Ongoing work is focussed on refining this approach with the aim of obtaining a better understanding of its advantages and limitations for future applications.






\section*{Appendix}
We report here all the parameters that were used for the numerical simulations reported in the paper.

The circular goal region is centered at the origin, i.e. $x_g=0$, with radius $\rho_g = 5\si{\meter}$.
The targets parameters were set to $\beta_1 = 1$, $\rho_\tau = 3\si{\meter}$, $v_{\tau}^{\max} = 9\si{\meter\second^{-1}}$. 
The random diffusive motion of the target uses $\Delta T = 1\si{\second}$ as update time interval and $\beta_{\max} = 1.8\si{\meter\second^{-1}}$ as maximum speed. 
The herder's maximum speed was set to $v_{h}^{\max} = 14\si{\meter\second^{-1}}$. The estimated radius of the influence zone assumed for the design of the control tutor is $\hat \rho_\tau = 1\si{\meter}$.
Each learning trial lasted $T_{tr}=100\si{\second}$ in the case of the single target experiments and $T_{tr}=500\si{\second}$ in the case of multiple targets, the sampling time was set to $T_s = 10^{-3} \si{\second}$.
The parameters of the learning update rule were set to  $\alpha = 0.9$ and  $\gamma = 0.8$ while the randomness parameter in the policies $\pi_Q$ and $\pi_T$ was set to $\varepsilon = 0.03$. 
The parameters of the reward function were set to $k_1 = 1$, $k_2 = 0.5$ and $\bar{k} = 100$.
To implement QL and CTQL, two alternative discretization of state space were tested. To reduce computational burden in the implementation the set of discretized relative distances were used to address and construct the Q-table.
A {\em coarser} discretization was obtained by sampling the range $[0,\hat \rho_\tau]$ of relative distances with stepsize $T_{m,d} := \frac{\hat \rho_\tau}{6}\si{\meter}$, and the range of angles [0,$2\pi$] with $T_{a,d} := {\frac{2\pi}{6}}\si{\radian}$.
The angular position of the herder was discretized in the range [0,$\frac{\pi}{2}$] with stepsize $T_{a,h} = \frac{\pi}{10}\si{\radian}$.
The target speed was discretized in the range $[0,v_h^{\max}]$ with stepsize $T_{m,v_\tau} = \frac{v_h^{\max}}{3}\si{\meter\second^{-1}}$, and the range of angles [0,$2\pi$] with $T_{a,v_\tau} = \frac{2\pi}{4}\si{\radian}$.
A finer discretization was obtained by reducing the sampling stepsizes of the quantities above to $T_{m,d} = \frac{\hat \rho_\tau}{10}\si{\meter}$, $T_{a,d} = {\frac{2\pi}{10}}\si{\radian}$, $T_{a,h} = \frac{\pi}{10}\si{\radian}$, $T_{m,v_\tau} = \frac{v_h^{\max}}{50}\si{\meter\second^{-1}}$, and $T_{a,v_\tau} = \frac{2\pi}{20}\si{\radian}$.
The action space consisted of possible herder velocities discretized in the range {[0,$v_{h}^{\max}$] with stepsize $T_{m,v_h} = \frac{v_{h}^{\max}}{10}\si{\meter\second^{-1}}$} and possible angular orientation in the range [0, $2\pi$] with $T_{a,v_h} = \frac{2\pi}{20}\si{\radian}$.
In the model of the target dynamics used for the control synthesis the parameter $\delta$ was set to unity while the tutoring control law gains were set to $k_i = 2$, $k_p =0.1$.

\end{document}